\documentclass[12pt]{article}
\usepackage{amsmath,amsfonts}
\usepackage{amsmath,amsfonts,latexsym}
\usepackage{amscd,amssymb,amsmath}
\usepackage{graphicx}
\newcommand{\ds}{\displaystyle}
\begin{document}
 \newtheorem{lema}{Lemma}[section]
\newtheorem{teo}[lema]{Theorem}
\newtheorem{propo}[lema]{Proposition}
\newtheorem{rema}[lema]{Remark}
\newtheorem{coro}[lema]{Corollary}
\newtheorem{defi}[lema]{Definition}
\newtheorem{exem}[lema]{Example}
\newtheorem{prop}[lema]{Proposition}
\newtheorem{assum}[lema]{Assumptions}
\newtheorem{abc}[lema]{}

\newcommand{\bea}{\begin{eqnarray*}}
\newcommand{\eea}{\end{eqnarray*}}

\newcommand{\pr}{\hspace*{6mm}}
\newcommand{\nd}{\noindent}

\newcommand{\af}{\alpha}
\newcommand{\ta}{\theta}

\newcommand{\hbx}{\hfill$\Box$}
\newcommand{\vv}{\vspace{4mm}}
\newcommand{\calc}{{\cal C}}
\newcommand{\Fix}{{\rm Fix}}
\newcommand{\reg}{{\rm reg}}
\newcommand{\ind}{{\rm ind}}

\title{When a smooth self-map of a  semi-simple Lie group  can   realize the least number of periodic points \\}
\author{Jerzy Jezierski}

 \maketitle

\centerline {Institute of Applications of  Informatics and Mathematics,}

 \centerline {Warsaw University of Life Sciences (SGGW)}

\centerline{Nowoursynowska 159, 00-757 Warsaw, Poland}

\centerline{ e-mail: $\; jerzy{\_}jezierski@sggw.pl$}

 \centerline{ AMS Subject Classification: Primary  55M20; Secondary  37C05, 37C25}

\centerline{ Keywords: Periodic points, Nielsen number, fixed point index, smooth maps, Lie group.}

\centerline{Research supported by the National Science Center, Poland,
  UMO-2014/15/B/ST1/01710.\\}

\begin{abstract}  There are two algebraic lower bounds of the number of $n$-periodic points of a self-map $f:M\to M$  of a compact smooth manifold of dimension at least $3$  :
$NF_n(f)=\min \{\# {\rm Fix}(g^n) ;g\sim f;\; g \text{ continuous}\}$ and $
NJD_n(f)=\min \{\# {\rm Fix}(g^n) ;g\sim f;\; g \text{ smooth}\}$. In general $NJD_n(f)$ may be much greater than $NF_n(f)$.
 We show that for a self-map of a semi-simple Lie group, inducing the identity fundamental group homomorphism, the equality
 $NF_n(f)=NJD_n(f)$
holds for all $n$ $\iff$ all eigenvalues of a quotient cohomology  homomorphism induced by $f$  have moduli $\le 1$.

\end{abstract}

\section{Introduction}
We consider a map $f:M \to M$, where $M$ is a compact connected manifold, ${\rm dim}M\ge 3$. We fix a number $n\in \mathbb{N}$ and we ask about the least number of n-periodic points
$$ \min_{g\sim f} \#{\rm Fix}(g^n)=?$$
in the homotopy class of $f$.

It turns out that the above minimum depends whether we consider   continuous or only  smooth maps homotopic to the given $f$. There are two algebraic homotopy invariants, for ${\rm dim} M\ge 3$,
 $$NF_n(f)=\min \{\# {\rm Fix}(g^n) ;g\sim f;\; g \text{ continuous}\}$$
 introduced by Boju Jiang in \cite{BJLectures} and
$$NJD_n(f)=\min \{\# {\rm Fix}(g^n) ;g\sim f;\; g \text{ smooth}\}$$
given in \cite{GJNonSimply}  . In general there may be a large gap between the two numbers $NF_n(f) \ll NJD_n(f)$. Even in the simply-connected case: $NF_n(f)=0 \text{ or } 1$ while $NJD_n(f)$ may be arbitrarily large \cite{SS}, \cite{GJsimply}.

 Nevertheless Cheng Ye You proved  that each self-map of a torus can be smoothly deformed to a map minimizing the number of $n$-periodic points which gives, in this case, the equality of two above invariants $NF_n(f) = NJD_n(f)$ \cite{You}. See also \cite{JeYou}.

It is natural to ask if such equality also holds for self-maps of other Lie groups. In \cite{Sinica} it was shown that each nonabelian Lie group admits a self-map  satisfying $NF_n(f) < NJD_n(f)$, proving that tori are the only Lie groups so that the equality holds for each self-map. Nevertheless we ask, which  self-maps of a Lie group $f:X\to X$ have the property:
{\it for each natural number $k$, there is a smooth map  $h$  homotopic to $f$ and  realizing the least number of $k$-periodic points?} In other words for which self-maps of Lie groups the equality $NF_k(f) = NJD_k(f)$ holds for all $k\in \mathbb{N}$.

It was shown in \cite{Sinica} that for a  self-map of a  compact Lie group with free fundamental group the above equality holds (for all $n$) if and only if the eigen-values of a cohomology homomorphism, defined by H.Duan  \cite{Duan}, have modulus $\le 1$. In this paper we show that this is also true for  self-maps of all compact connected semi-simple Lie groups inducing the identity isomorphism of fundamental groups.

In Section \ref{SectionIndices} we give a necessary background about the sequences of indices of iterations of a map. In Section \ref{SectionNielsen} we outline the Nielsen fixed point theory. In Section  \ref{SectionDuan} we recall the class of rational exterior spaces which contains Lie groups and is very convenient here. In the last Section we prove the main Theorem \ref{MainTeoIdentity}. Here is the main idea of the paper: Theorem \ref{CritForSmooth=Cont} reduces the problem of minimizing the number of periodic points by a smooth map to finding a system of expressions attached to some vertices in the Reidemeister Graph and realizing the index function ${\rm ind}(f^b;B)$ for each orbit of Reidemeister classes. On the other hand Lemma  \ref{L(f^)SRImpliesGraph} guaranties a system realizing the index function for a part of the Graph. In the last Section we show that under our assumptions the above system also realizes the index function for all vertices of the Reidemeister graph.

\begin{rema}
 If $f:M\to M$ is a smooth map then it is not important what homotopies (continuous of smooth) between $f$ and $g$ we use in the definition of $NJD_n(f)=\min \{\# {\rm Fix}(g^n) ;g\sim f;\; g \text{ smooth}\}$. In fact any continuous homotopy between two snooth maps can be replaced by a smooth homotopy: first we make the homotopy constant for $0\le t\le \epsilon$ and then for $1-\epsilon\le t\le 1$ and finally we corrrect the homotopy inside $M\times (0,1)$.

\end{rema}



\section{Indices of iterations of a smooth map}\label{SectionIndices}

   Albrecht Dold \cite{Dold}  noticed that a sequence of fixed point
  indices $A_k={\rm ind}(f^k;x_0)$, where $f$ is a continuous self-map of a
  Euclidean space $ \mathbb{R}^m$ and $x_0$ is an isolated fixed point for each
  $f^k$, must satisfy some congruences. Namely for each $n\in \mathbb{N}$

  $$\sum_{k|n} \mu(n/k)\cdot {\rm ind}(f^k;x_0) \equiv 0 (modulo\; n)$$
  where $\mu$ denotes the M$\ddot{o}$bius function.

  It was shown \cite{BaBo} that each sequence of integers $(A_k)$ satisfying Dold
    congruences can be realized as $A_k={\rm ind}(f^k;x_0)$, for a continuous
    self-map of $\mathbb{R}^m$ for  $m\ge 3$. In other words Dold congruences are
     the only restrictions for the sequence of the fixed point index of a
     continuous map.

  Surprisingly it turned out that there are much more restrictions on  sequences
   $A_k={\rm ind}(f^k;x_0)$ when $f$ is smooth  \cite{SS}, \cite{CMPY}.

  \begin{defi}
 A sequence of integers $A_n$  will be called {\it smoothly realizable in}
 $\mathbb{R}^m$ (or {\it in dimension m}) if there exists a smooth self-map
 $f: \mathbb{R}^m \to \mathbb{R}^m$ and an isolated fixed point
 $x_0\in {\rm Fix}(f)$, which is also an isolated fixed point of each iteration
  $f^n$, so that $A_k={\rm ind}(f^k;x_0)$.

  \end{defi}

    In  Theorem \ref{smooth.real.teo} we recall all possible sequences which can be
    obtained as fixed point indices of  a smooth self-map of $\mathbb{R}^m$
    (for $m\ge 3$).

  It is convenient to present the sequences of integers as the sum of the
  following elementary periodic sequences
  \begin{defi} For a given $k \in {\mathbb N}$ we define
$$\reg_k(n)= \begin{cases} k &  \text{if} \;\; k \mid n, \\
0 & \text{if}\;\; k \nmid n.
\end{cases} $$
\end{defi}

In other words, $\reg_k$ is the periodic sequence:
$$(0, \ldots, 0, k, 0, \ldots, 0,k, \ldots),$$ where the
non-zero entries appear for indices divisible by $k$.

It turns out that
each integer sequence  $(A_n)$ can be written down
uniquely in the following form of a periodic expansion:
$ A_n= \sum_{k=1}^{\infty} a_k \reg_k(n),$ where
$a_n=\frac{1}{n} \sum_{k|n} \mu(\frac{n}{k})\; A_k$.
 Moreover all  coefficients $a_k$ are integers if and only if the sequence
 $(A_n)$ satisfies Dold congruences.

 \begin{rema}
 Later we will be interested in iterations $f^k$, for $k$ dividing a prescribed
 $n$, so we will consider only finite sequences labeled by $\{k\in \mathbb{N}; k|n\}$.
 We will say that a finite sequence $(B_k)_{k|n}$ is smoothly realizable in
 $ \mathbb{R}^m$ if it is the restriction of a  sequence smoothly realizable
   in $\mathbb{R}^m$. \hbx

 \end{rema}

We will use the following notation. For a finite subset $A\subset \mathbb{N}$ we denote by ${\rm lcm}(A)$ the least common multiplicity. Moreover we will define
${\rm LCM}(A)=\{{\rm lcm}(B); B\subset A\}$ and  ${\rm LCM}_2(A)={\rm LCM}(A\cup\{2\})$.
\begin{teo}\label{smooth.real.teo} Thm. 2.5 in \cite{{JeMod_2}}\\
A sequence $(D_n)$ is smoothly realizable in dimension $m$ if and only if there exist natural numbers $d_1,...,d_s$ ($d_i\ge 3$ ,$2s\le m$)  so that
\begin{equation}
\label{FormulaInMainThm}
    D_n= \sum_{k} \alpha_k\reg_k(n)
\end{equation}
where the summation runs through the set  $\text{ LCM }_2(\{d_1,...,d_s\}$ and the coefficients $\alpha_k$ are integers.
If moreover  $m\le 2s+2$ then   the following restrictions hold
\begin{enumerate}
  \item[[0]] if $m=2s$ then $\alpha_1=1$ and $LCM$. Here $LCM$ is the shorthand for : the restriction  of summation (1)  must run through
$\text{ LCM }(\{d_1,...,d_s\})$ (without this restriction the summation may run through
$\text{ LCM }_2(\{d_1,...,d_s\})$).
  \item[[1]] if $m=2s+1$ then ($|\alpha_1|\le 1$  and   $LCM$) or ($\alpha_1=1$ and ($\alpha_2=0$  or $1$))
  \item[[2]] if $m=2s+2$ then  $|\alpha_1|\le 1$ or  $LCM$.
\end{enumerate} \hbx

\end{teo}

\begin{coro}\label{smooth.real.Coro}
If the above sequence $(D_n)$ is smoothly realizable in dimension $m$ then
\begin{enumerate}
  \item $|\alpha_1|\ge 2$ implies $m\ge 2s+3$ or ($m=2s+2$ and $LCM$).
  \item $\alpha_1=0$ implies  $m\ge 2s+2$ or ($m\ge 2s+1$ and $LCM$).
\end{enumerate}\hbx

\end{coro}



  It is easy to notice that if  $f:M\to M$ is a self-map of a compact manifold  and ${\rm Fix(f^n)}$ is a point then the sequence $(L(f^k))_{k|n})$ is smoothly realizable in  $\mathbb{R}^m$ where $m = {\rm dim}M$. It turns out that if $M$ is simply-connected and  ${\rm dim}(M) \ge 3$ then the inverse implication is also true.\cite{GJsimply}

\section{Nielsen fixed point theory}\label{SectionNielsen}

We consider a self-map of a compact connected polyhedron $f:X\to X$ and its fixed point set ${\rm Fix}(f)$. We define the {\it Nielsen relation} on this set by :

$x\sim y$ if there is a path  $\omega$ joining $ x $ with  $y$  so that   $f\omega$ and  $\omega$ are fixed end point homotopic\\
This relation splits ${\rm Fix}(f)$ into {\it Nielsen  classes}. Their set will be denoted by $\mathcal{N}(f)$. We say that a Nielsen class $A$ is ${\it essential}$ if its fixed point index is nonzero : ${\rm ind}(f;A)\ne 0$. The number of essential Nielsen classes is called {\it Nielsen number} and denoted $N(f)$. This is a homotopy invariant and moreover it is the lower bound of the number of fixed points in the (continuous) homotopy class: $N(f)\le \min_{h\sim f} \# {\rm Fix}(h)$ \cite{Brown}, \cite{BJLectures}, \cite{JMBook}.

 On the other hand we define, the set of {\it Reidemeister classes} of the map $f$ as the quotient set of the action of the fundamental group $\pi_1M$
on itself given by  $\omega * \alpha=  \omega\cdot \alpha \cdot (f_{\#}\omega)^{-1}$. Here we take as the base point a fixed point of $f$. We denote the quotient space by $\mathcal{R}(f)$. There is a natural injection from the set of the Nielsen classes to the set of Reidemeister classes $ \mathcal{N}(f)\subset \mathcal{R}(f)$ defined as follows. We choose a point $x$ in the given Nielsen class $A$ and a path $\omega$ from the base point $x_0$ to $x$. Then the loop $\omega*(f\omega)^{-1}$ represents the corresponding Reidemeister class.

Now we consider the iterations of the map $f$. For fixed natural numbers $l|k$ there is a natural inclusion ${\rm Fix}(f^l)\subset  {\rm Fix}(f^k)$ which induces the map $\mathcal{N}(f^l)\to \mathcal{N}(f^k)$ (which may be not injective). This map  extends to $i_{kl}: \mathcal{R}(f^l) \to \mathcal{R}(f^k)$ and the last is given by
$$i_{kl}[x]=[x\cdot f^l(x)\cdot f^{2l}(x)\cdots f^{k-l}(x)]$$
The functorial equalities are satisfied: $i_{kl}i_{lm}=i_{km}$ , $i_{kk}={\rm id}$ and moreover the diagram
$$ \begin{CD}
  \mathcal{N}(f^l) @> i_{kl} >> \mathcal{N}(f^k)\\
  @V  VV @VV   V\\
   \mathcal{R}(f^l) @> i_{kl} >> \mathcal{R}(f^k)\\
 \end{CD}
 $$
commutes.

 The group $\mathbb{Z}_k$ acts on ${\rm Fix}(f^k)$ by
$$ {\rm Fix}(f^k)\ni [x] \to [fx]\in    {\rm Fix}(f^k)$$
and on $\mathcal{R}(f^k)$ by
$$ \mathcal{R}(f^k) \ni [a]\to [f_{\#}(a)]\in    \mathcal{R}(f^k)$$
Then the diagram
$$ \begin{CD}
  {\rm Fix}(f^k) @ >>> {\rm Fix}(f^k)\\
  @V  VV @VV   V\\
   \mathcal{R}(f^k) @>>> \mathcal{R}(f^k)\\
 \end{CD}
 $$ commutes.
  We denote by $ \mathcal{OR}(f^k)$  the set of orbits of the above action {\it (orbits of Reidemeister classes)}.

 Consider Reidemeister classes $A\in \mathcal{R}(f^k)$, $B\in \mathcal{R}(f^l)$ , $l|k$ , satisfying  $i_{kl}(B)=A$. Then we say that $A$ {\it reduces} to $B$ , or $B$ {\it preceeds} $A$ and we write $B \preceq A$. The class is called {\it reducible} if $B \preceq A$ for a $B\ne A$. A similar definition works  for orbits.

Let $\mathcal{ER}(f^k)$ , $\mathcal{IR}(f^k)$ denote sets of essential and irreducible classes, respectively.

A map $f$ is called {\it essentially reducible} if each class $A\in \mathcal{R}(f^l)$   preceding an essential class  is essential.

We will denote by $\mathcal{IEOR}(f)$, or simply $\mathcal{IEOR}$, the set of irreducible essential orbits of Reidemeister classes of $f$.

\begin{coro}
   Lie groups  are essentially reducible \cite{BJLectures}.
\end{coro}

 It is convenient to put all Reidemeister data into a directed graph $\Gamma$.

\begin{itemize}
  \item Vertices are elements of the  union ${\rm Vert}(\Gamma)= \bigcup_{k\in \mathbb{N}} \mathcal{OR}(f^k)$.
  \item There is a unique directed edge from $A\in \mathcal{OR}(f^l)$ to $B\in \mathcal{OR}(f^k)$ if $i_{kl}(A)=B$.

  \end{itemize}

  Moreover to each vertex $A\in \mathcal{OR}(f^k)$ an integer ${\rm ind}(f^k;A)$ is defined, and the following {\it Dold Congruences} are satisfied

\begin{lema} (Lemma 3.3 in \cite{GJNonSimply})

   For each $B\in \mathcal{OR}(f^b)$  :
  $$\sum_{C \preceq B} \mu(b/c)\cdot {\rm ind}(f^c:C) \equiv 0 (modulo\; b)$$
  where $\mu$ denotes the Mobius function and $C\in \mathcal{OR}(f^c)$. \hbx
  \end{lema}

\begin{defi}\label{DefOfReg}
 \begin{enumerate}
   \item For a fixed  $A\in \mathcal{OR}(f^k)$ and a number $r\in \mathbb{N}$ we define
${\rm Reg}_A^r:{\rm Vert}(f) \to \mathbb{Z}$ by the formula
 $${\rm Reg}_A^r(B)=   \begin{cases}
k\cdot r & \text{ if } i_{rk,k}(A) \preceq B \\
0 & \text{ otherwise }
\end{cases} $$

  \item More generally for each Dold sequence $d_k =\sum_l a_l{\rm reg}_l(k)$ and an orbit $A\in \mathcal{OR}(f^r)$ we define a function $C:{\rm Vert}(f^k)\to \mathbb{Z}$ as $C(B) =\sum_l a_l{\rm Reg}_A^l(B)$. We say then that the sequence $\sum a_l{\rm reg}_l$ is {\it attached} to the orbit $A$ or that $D: {\rm Vert}(f^k)\to \mathbb{Z}$ {\it comes from } $\sum a_l{\rm reg}_l$.
\end{enumerate}
\end{defi}

\begin{coro}\label{LeastCont}(Corollary 7.2 in \cite{Sinica}).
A continuous  essentially reducible map $f:M\to M$ realizes the least number of $n-$periodic points $\iff$ for each orbit  $A\in \mathcal{IEOR}(f^k)$ ($k|n$), the corresponding orbit of Nielsen classes (also denoted by $A$)  contains exactly one $k$-orbit of points ($A=\{a_1,...,a_k\}$) and all other orbits are empty. \hbx
\end{coro}

Now we are in a position to make precise when the map in the above Corollary may be smooth.

\begin{teo}\label{CritForSmooth=Cont}(Theorem (7.5) in \cite{Sinica}).
  An essentially reducible map $f:M\to M$ is homotopic to a smooth map $g$ realizing the least number of $n-$periodic points $\iff$ one can attach to  each orbit $A\in \mathcal{IEOR}(f^k)$ ($k|n$) an expression $C_A$, realizable in $\mathbb{R}^m$,  so that
  $${\rm ind}(f^k;B)= \sum_{A\in \mathcal{IEOR}} C_A(B)$$
  for each $B\in \mathcal{OR}(f^k)$  , $k|n$. Here $C_A: {\rm Vert}(f)\to \mathbb{Z}$ is the function given by Definition \ref{DefOfReg}\hbx
\end{teo}

  The above Theorem and Theorem $ \ref{smooth.real.teo}$ allow   to determine whether  the least number of $n-$periodic point can be realized by a smooth map.

 \begin{defi}
 If the right hand side in Theorem $\ref{CritForSmooth=Cont}$ holds then we say that the Reidemeister graph $\mathcal{GOR}(f)$  is {\it smoothly realizable in dimension $dim(M)$ for ${k|n}$}. By  Theorem $\ref{CritForSmooth=Cont}$  this is equivalent to a smooth minimization of ${\it Fix}(f^n)$

 \end{defi}

  We end the Section by a Remark which guaranties the partial realization of the index function.

  \begin{rema}\label{L(f^)SRImpliesGraph}
  We fix a map $f:M\to M$ and a number $d\in \mathbb{N}$ so that each sum $\sum_{l|d}a_l\cdot l$ is smoothly realizable in dimension $m={\rm dim}M$ for $k|d$. We moreover assume that $f$ is essentially reducible for $k|d$. Then the Reidemeister graph $\mathcal{GOR}(f)$  is smoothly realizable for $k|d$.

  \end{rema}
  {\it Proof}. This is intuitively clear. Since each essential orbit reduces to an orbit in $\mathcal{IEOR}$, we may subordinate (in any way) to each essential orbit an orbit in  $\mathcal{IEOR}$ preceding it. This defines expressions  attached in $\mathcal{IEOR}$ realizing all essential classes. Here we consider only orbits for $k|d$, so the assumptions on the sums implies that the obtained expressions are smoothly realizable in $\mathbb{R}^m$.

      Now we present a formal proof.

   We need to attach to  each orbit $A\in \mathcal{IEOR}(f^k)$ ($k|n$) a sequence $C_A$, realizable in $\mathbb{R}^m$ and satisfying the equality in  Theorem \ref{CritForSmooth=Cont}.

   Let $Vert'$ denote the set of  Reidemeister orbits reducing to an essential class. Of course each essential orbit belongs to $Vert'$ and by essential reducibility each orbit in $Vert'$ reduces to an orbit in $\mathcal{IEOR}$.
    For each $B\in \mathcal{OR}(f^b)\cap Vert'$ we fix an $A_B \in \mathcal{IEOR}(f^{a_b})$ satisfying $A_B\preceq B$. Then for any $E\in \mathcal{OR}(f^e)$.

  $${\rm ind}(f^e;E)=\sum_{B\in Vert'} a_B\cdot {\rm Reg}_B^1(E)= \sum_{A\in \mathcal{IEOR} }\left(\sum_{B\in Vert'; A_B=A} a_B\cdot {\rm Reg}_B^1(E)\right)$$

  Now it remains to show that for each $ A\in \mathcal{IEOR}$ the expression $\ds C_A= \sum_{B\in Vert'; A_B=A} a_B\cdot {\rm Reg}_B^1$ comes from a sequence smoothly realizable in $\mathbb{R}^m$ for $c|d$.

  We notice that $\ds \sum_{B; A_B=A} a_B\cdot {\rm Reg}_B^1$ is a subsum of $\ds \sum_{B; A\preceq B} a_B\cdot {\rm Reg}_B^1$ so it remains to prove the claim for the last sum. We notice that for a given orbit $A$
    $$\{B; A\preceq B\} = \{i_{la,a}(A); l\in \mathbb{N}\}$$

    Now $\ds \sum_{B;A\preceq B} a_B\cdot {\rm Reg}_B^1$ comes from
     the expression $\ds \sum_{l=1}^{d/a} a(l){\rm reg}_l$ attached at $A$,
     where $a(l)=a_B$ for $B=i_{la,a}(A)$. \hbx

  \section{Rational exterior spaces}\label{SectionDuan}
  Let us recall recall that the Lefschetz number of a self-map of a compact Lie group is given by $L(f^k)={\rm det}(I-A^k)$ for an integer matrix $A$. In fact such equality for a larger class of spaces introduced by Haibao Duan  \cite{Duan}. Here is a short description of the matrix $A$.

  We are given a topological space $X$ and we consider the rational cohomology $H^*(X;\mathbb{Q})$. An $x\in H^r(X;\mathbb{Q})$, is called {\it decomposable} if $x=\sum_i x_i\cup y_i$ for some $x_i\in H^{p_i}(X;\mathbb{Q})$ , $y_i\in H^{q_i}(X;\mathbb{Q})$ , $p_i+q_i=r$ , $p_i\ge 1, q_i\ge 1$. Let $D^r(X)\subset H^r(X;\mathbb{Q})$ denote the space over $\mathbb{Q}$ consisting of all decomposable elements. Then the quotient $A^r(X)=H^r(X;\mathbb{Q})/D^r(X)$ is a vector space over $\mathbb{Q}$. A continuous map $f:X\to Y$ induces $f^*_r: H^r(Y; \mathbb{Q}) \to H^r(X; \mathbb{Q})$ and $A^r(f):A^r(Y) \to A^r(X)$ for each $r\ge 0$. Let $A(X)=\oplus_{r=1}^{\infty}A^r(X)$ and $A(f):A(Y)\to A(X)$ be the induced homomorphism.

  \begin{exem}
  If $H^*(M;\mathbb{Q})=\Lambda(a_{i_1},...,a_{i_s};\mathbb{Q})$ is an exterior algebra then $A^*(M)= \mathbb{Q}(a_{i_1},...,a_{i_s})$ is a free $\mathbb{Q}$ linear space and $A^r(M)$ is the subspace spanned by these above generators which satisfy $a_i\in H^r(M; \mathbb{Q})$.

  \end{exem}

  \begin{defi}\cite{Duan}
  A connected topological space $X$ is called {\rm rational exterior} if some homogeneous elements $\omega_i\in H^{odd}(X; \mathbb{Q})$ can be chosen so that the inclusions $\omega_i\to H^*(X; \mathbb{Q})$ give rise to a ring isomorphism $\Lambda_{\mathbb{Q}}(\omega_1,...,\omega_k)= H^*(X;\mathbb{Q})$. \hbx

  \end{defi}

  The examples of exterior rational spaces are given in \cite{Duan} pages 73-74. By a Hopf Theorem  compact Lie groups are exterior rational spaces \cite{Hopf}.

  \begin{rema}
  If $X$ is a rational exterior manifold and $H^*(X;\mathbb{Q})= \Lambda_{\mathbb{Q}}(\omega_1,...,\omega_n)$ then $\omega_1 \cup ... \cup \omega_n \ne 0$,  which implies ${\rm dim}(\omega_1)+...+{\rm dim}(\omega_n)\le {\rm dim}M$. \hbx

  \end{rema}
  \begin{teo}\label{Duan}\cite{Duan} (Duan)\\
  Let $f$ be a self-map of a rational exterior space. Then $L(f^n)= det(I-(A(f))^n)$ for all $n \in \mathbb{N}$. \hbx

  \end{teo}

\begin{rema}\label{Less=1Implies=1}
   Let $L(f^n)= det(I-(A(f))^n)$. We define $\bar A^k(X)= A^k(X)/{\rm ker}A^k(f)$ , $\ds \bar A(X)= \oplus_{k=1}^{\infty} \bar A^k(X)$ and
   $\ds \bar A(f)= \oplus_{k=1}^{\infty} \bar A^k(f)$ .
   Then $det(I-(A(\bar f))^n)=det(I-(A(f))^n)=L(f^n)$. Thus we may assume that the matrix $A(f)$ in the equality $L(f^n)= det(I-(A(f))^n)$ may have only eigenvalues of modulus equal to $1$.

\end{rema}

\section{Finite fundamental group , $f_{\#}=identity$ }\label{SectionFinal}

In this Section we will prove the main result of the paper. Let us recall that each   Lie group with the finite fundamental group $\pi_1M$ is at least three dimensional, so  algebraic invariants ${\rm NF}_n(f)$ and ${\rm NJD}_n(f)$ are respectively the least lower bounds of the number of periodic points in the continuous  and smooth homotopy class of $f$.

Let $q_1<\cdots < q_s$ be the primes dividing the order of the finite abelian group $\pi_1M$. Then

\begin{equation}\label{DecompositionOfpi_1(M)}
  \pi_1M= G(q_1)\oplus \cdots \oplus G(q_s)
\end{equation}
 where $G(q) = \mathbb{Z}_{q^{\alpha_1}} \oplus \cdots \oplus \mathbb{Z}_{q^{\alpha_u}}$ , $\alpha_i\in \mathbb{N}$.

 \begin{teo}\label{MainTeoIdentity}  We are given a self-map of a compact connected semi-simple Lie group $f: M\to M$   inducing the identity homomorphism
 of the fundamental group. Then  the least number of n-periodic points in the homotopy class of $f$ can be realized by a smooth map $g_n$ (for each  a priori fixed number $n\in \mathbb{N}$) if and only if either $L(f)=0$ or all eigenvalues of $A(f)$ have moduli $\le 1$.
  \end{teo}
  {\it Proof}.   $ \Longrightarrow$
  We assume that $L(f)\ne 0$ and an eigen-value of $A(f)$ has modulus $>1$. We will show that then $\mathcal{GOR}(f^n)$ is not smoothly realizable  for some $n\in \mathbb{N}$ contradicting to the assumption.

   Let $q_1<\cdots <q_s$ be the prime numbers involved in the decomposition of $\ds \pi_1M$: see (\ref{DecompositionOfpi_1(M)}).

   On the other hand we denote by $d_1,...,d_s$
    the minimal periods of roots of unity in spectrum of $A(f)$.


  Since $q_1,...,q_s$ are the only primes dividing $\# \pi_1 M$, by Lemma \ref{L(f^n)TendsToInfty} it is enough to find  a sequence $(r_k)$ of natural numbers satisfying:
      (1) each $r_k$ is relatively prime with $q_1,...,q_s$ and
        (2)  $\ds \lim_{k\to \infty}|L(f^{r_k})|=\infty$. If we assume that $\lim_{n\to \infty}|L(f^k)|=\infty$  then it is enough to take any sequence of numbers relatively prime with $q_1,...,q_s$. .
  In the general case, under the assumptions of the Theorem the sequence $|L(f^k)|$ need not to tend to infinity. To see this  we regroup the eigenvalues so that:

          $|\lambda_i| > 1$ for $1\le i<s_1$ ,

          $|\lambda_i|= 1$ for $s_1\le i<s_2$ and $\lambda_i$ is not a root of unity ,

           $|\lambda_i|= 1$ for $s_2\le i<s_3$ and $\lambda_i$ is  a root of unity ,

           $|\lambda_i|< 1$ for $s_3\le i<s$,

           for some $1\le s_1\le s_2\le s$. Here $s_1\ge 2$, since by the  assumption an eigen-value has modulus $>1$. Now

           \begin{equation}\label{lambda}
                        |L(f^k)|= |(1-\lambda_1^k)\cdots (1-\lambda_s^k)|=
               \end{equation}
             $$|\lambda_1^k \cdots   \lambda_{s_1-1}^k|
            \cdot |(1/(\lambda_1^k) -1) \cdots (1/(\lambda_{s_1-1}^k) -1)|
            \cdot |(1-\lambda_{s_1}^k) \cdots (1-\lambda_{s_2-1}^k)|$$
            $$ \cdot |(1-\lambda_{s_2}^k) \cdots (1-\lambda_{s_3-1}^k)|
            \cdot  |(1-\lambda_{s_3}^k) \cdots (1-\lambda_{s}^k)|$$

            We notice that in the last product the first factor tends to infinity (by the assumption that $s_1\ge 2$), the second and the last tend to $1$.
            If the third and the fourth factors were bounded from below by a positive number, the whole product would tend to infinity and the Theorem would be proved by the first part.
            However the fourth factor corresponds to eigenvalues whose some powers are equal $1$ so the corresponding Lefschetz number is zero. On the other hand $\lambda_s^k$ which occur in the third factor (for $s=s_1,...,s_2-1$), may  be very close to $1$  which can make the product very small. We will overcome this difficulty by finding a sequence $(r_k)$ of natural numbers  relatively prime with $q_1,...,q_u;d_1,...,d_s$ and such that the numbers $\lambda_{s_1}^{r_k}, \cdots , \lambda_{s_2-1}^{r_k}$ avoid an arc $C_{\epsilon}\subset S^1$ containing $1$, for all $k$. The last will make the third factor of (\ref{lambda}) bounded from below by a positive number. Then it remains to notice that $\lambda_{s_2}^{k},..., \lambda_{s_3-1}^{k}$ are cyclic so they reach only a finite number of values. Moreover no one of these powers equals 1, since chosen  numbers is are relatively prime   with all $d_1,...,d_s$. Let us emphasize that the choice of the sequence will be possible because the assumption $L(f)\ne 0$ guarantees that all $d_i\ne 1$. Since on the other hand  $q_j\ne 1$, there are infinitely many numbers relatively prime with all $q_1,...,q_u;d_1,...,d_s$.

            Now we show the existence of a sequence and numbers satisfying the above conditions.

   Let us denote $L:= lcm(q_1,...,q_u; d_1,...,d_s)$ and let
    $\phi(L)=\#\{k\in \mathbb{N}; k\le L , {\rm gcd}(k,L)=1\}$ be the Euler totient function. Then
  $$\lim_{n\to\infty}\frac{\#\{k\in \mathbb{N}; k\le n , {\rm gcd}(k,L)=1\}}{n}=\frac{\phi(L)}{L}>0$$
   Let us recall that if $\lambda = exp(2\pi\omega i)$ where $\omega$ is irrational then the sequence $\lambda^k$ is equidistributed \cite{Fourier} i.e. for each arc $C\subset S^1$
   $$\lim_{n\to \infty} \frac{\#\{k\in \mathbb{N}; \lambda^k\in C, k\le n\}}{n}= \frac{\text{length of } C}{{\text{length of } S^1}}$$

   Now it is enough to prove
     that there exists $\epsilon >0$ so that for infinitely many $k\in \mathbb{N}$  relatively prime with $L$ the following holds:

             $$\text{ no number}  \lambda^{k}_{s_1},..., \lambda^{k}_{s_2-1} \text{ belongs to the arc } C_{\epsilon}=\{{\rm exp}(2\pi\omega  i) ; |\omega|\le \epsilon\} $$

             We fix an $\epsilon>0$ and we notice that
             $$\lim_{n\to \infty} \frac{\#\{k\in \mathbb{N}; \lambda^k\in C_{\epsilon}, k\le n\}}{n}= \frac{\text{length of } C_{\epsilon}}{{\text{length of } S^1}}=\frac{2\pi\cdot 2\epsilon}{2\pi}=2\epsilon$$
             This implies

             $$\limsup_{n\to \infty} \frac{1}{n}\#\{k\in \mathbb{N};  k\le n, \lambda_i^k\in C_{\epsilon}  \text{ for an } l\in\{s_1,...,s_2-1\}\} $$

              $$\le \frac{1}{n}\left( \sum_{l=s_1}^{s_2-1} \limsup_{n\to \infty} \#\{k\in \mathbb{N}; k\le n , \lambda_l^k\in C_{\epsilon} \}\right)= (s_2-s_1)\cdot 2\cdot \epsilon$$

                  This implies
              $$\liminf_{n\to \infty} \frac{1}{n}\#\{k\in \mathbb{N};  k\le n, \lambda_l^k\notin  C_{\epsilon}  \text{ for all } l\in\{s_1,...,s_2-1\}\}\ge 1-2(s_2-s_1)\epsilon $$
              We get
              $$\liminf_{n\to \infty} \frac{1}{n}\#\{k\in \mathbb{N};  k\le n , \lambda_l^{k}\notin  C_{\epsilon}  \text{ for all } l\in\{s_1,...,s_2-1\} \text{ and } {\rm gcd}(k,L)=1 \}$$
              $$ \ge \left(1-(s_2-s_1)\epsilon +\frac{\phi(L)}{L} \right)-1=\frac{\phi(L)}{L}-2(s_2-s_1)\epsilon$$
              Now if we put $\epsilon$ sufficiently small, say $\epsilon=\frac{\phi(L)}{4L(s_2-s_1)}$, we get the last limit is positive which proves that the set of numbers $k\in \mathbb{N}$ satisfying: $\lambda_s^k\notin C_{\epsilon}$ (for  $s_1\le s<s_2$) and $k$ is relatively prime with $L$, is infinite.
               This gives the desired sequence $(r_k)$.

 $\Longleftarrow$
  Let  $d_1,...,d_s$ be the minimal periods of roots of unity in spectrum of $A(f)$ and let $d=lcm(d_1,...,d_s)$. Then the sequence $(L(f^k))_{k|d}$ is smoothly realizable in $\mathbb{R}^m$, see Remark \ref{L(f)IsSRealizable}.  Now  Lemma \ref{L(f^)SRImpliesGraph} implies  the existence of expressions $\{C_A\}_{A\in \mathcal{IER}(f^a)}$ smoothly realizable in $\mathbb{R}^m$ satisfying
    \begin{equation}\label{DoldDecompForId}
           {\rm ind}(f^b;B)= \sum_{A\in \mathcal{IEOR}} C_A(B)
     \end{equation}
 for  $B\in \mathcal{R}(f^b)$ with $b|d$. Let us notice that the above equality also holds for $b|d$  if in the sums $C_A=\sum_{E\in \mathcal{IEOR}}a_E {\rm Reg}_E^1$ we drop all summands $a_E {\rm Reg}_E^1$  , $E\in \mathcal{R}(f^e)$,  $e\nmid d$. In other words we may assume that $a_E\ne 0$ implies $e|d$.

    To prove the Theorem it is enough to show that the equality  (\ref{DoldDecompForId}) holds for all $B$. We consider two cases.

 1. First we assume that for $\bar b=gcd(b,d)$ the map $i_{b,\bar b}$ is the bijection.
 Let $\bar B\in \mathcal{R}(f^{\bar b})$ be the unique class preceding $B$.
 By Lemma \ref{L=L^bar} $L(f^b)=L(f^{\bar b})$, hence Jiang implies ${\rm ind}(f^b,B)={\rm ind}(f^{\bar b},\bar B)$. It remains to show that
 $\sum_{A\in \mathcal{IEOR}}C_A(B)= \sum_{A\in \mathcal{IEOR}}C_A(\bar B)$. Since $\sum_A C_A$ is the sum of ${\rm Reg}_E^1$ where $e|d$, it remains
  to notice that ${\rm Reg}_E^1(B)={\rm Reg}_E^1(\bar B)$. The last holds for $e|d$, since then $E\preceq B \iff E\preceq \bar B$.

  2.  Now we assume that $i_{b,\bar b}$ is not a bijection. By Lemma \ref{bar_kIZO} this occurs only for $b$ inessential. Then $0=L(f^b)=L(f^{\bar b})$,
  hence Jiang implies ${\rm ind}(f^b;B)=0$ and ${\rm ind}(f^{\bar b};\bar B)=0$ for all $\bar B\in \mathcal{OR}(f^{\bar b})$. Finally it remains to show  that
   $\sum_{A\in \mathcal{IEOR}}C_A(B)=0$. In fact
   $$\sum_{A\in \mathcal{IEOR}}C_A(B)= \sum_{\bar B\in \mathcal{R}(\bar f^b)}(\sum_{A\in \mathcal{IEOR}; A\preceq \bar B}C_A(B))=
   \sum_{\bar B\in \mathcal{R}(\bar f^b)} {\rm ind}(f^{\bar b}) =\sum 0=0$$ \hbx



      It remains to prove Lemmas \ref{L(f^n)TendsToInfty}, \ref{bar_kIZO} ,\ref{L=L^bar}.

      \begin{lema}\label{L(f^n)TendsToInfty}
      We are given a self-map $f: M\to M$,  of a compact connected  Lie group with $\pi_1M$ finite,  inducing the identity homomorphism
 of the fundamental group. If moreover there exists a sequence $(r_k)$ of natural numbers satisfying
      \begin{enumerate}
        \item each $r_k$ is relatively prime with $\#\pi_1M$
        \item $\ds \lim_{k\to \infty}|L(f^{r_k})|=\infty$
      \end{enumerate}
      then there exists a natural number $n$ so that $\mathcal{GOR}(f^k)_{k|n}$ is not smoothly realizable in dimension ${\rm dim}(M)$.

      \end{lema}
      {\it Proof}; Let $q_1,...,q_s$ be the primes dividing  $\#\pi_1M$.
   Since $f_{\#}=id$, the maps $i_{kl}:\mathcal{R}(f^l)\to \mathcal{R}(f^k)$ are given by
   $$i_{kl}[x] =  (id+f^l_{\#}+f^{2\cdot l}_{\#}+\cdots + +f^{k-l}_{\#})[x]=[x+\cdots+x]=[k/l\cdot x]$$
 In particular if $k/l$ is relatively prime with all $q_1,...,q_u$ then $i_{kl}$ is a bijection. Now for each number $r$ relatively prime with $q_1,...,q_u$ the Reidemeister class (=orbit) $[0]\in \mathcal{R}(f^1)$ is the unique essential irreducible class to which reduces $[0]\in \mathcal{R}(f^r)$. Let us fix a number $v\in \mathbb{N}$ and let  $r_0=lcm(r_1, \cdots, r_v$). If $(\mathcal{GOR}(f))_{k|r_0}$ were smoothly realizable in $\mathbb{R}^m$, then one could attach to the class $[0]\in \mathcal{R}(f^1)$ a smoothly realizable expression realizing all
   ${\rm ind}(f^{r_i};[0])$ for $i=1,...,v$. Since by the assumption $|L(f^{r_i})|$ tends to the infinity, ${\rm ind}(f^{r_i};[0])=\frac{1}{\#\pi_1M}|L(f^{r_i}|$ also tends to infinity. Now taking sufficiently large $v$ we may get arbitrarily many different values in the sequence ${\rm ind}(f^{r_i};[0])$ for  $i=1,...,v$.
   But any smoothly realizable expression in $\mathbb{R}^m$ may realize at most $2^{[(m+1)/2]}$ values \cite{BaBo}. The obtained contradiction means that $(\mathcal{GOR}(f^k))_{k|r_0}$ can not be smoothly realized for sufficiently large $r_0$. \hbx

      \begin{rema}\label{L(f)IsSRealizable}
       Let $f:M\to M$ be a self-map of a compact connected   Lie group  with $\pi_1M$ finite and such that all eigenvalue have moduli $\le 1$. Let $d_1,...,d_s$ be  minimal degrees of eigen-values of $A(f)$ and let $d= lcm(d_1,...,d_s)$. By \cite{Sinica} (page 1482, points (1),(2),(3)) $(L(f^k)_{k|d})$ is smoothly realizable in $\mathbb{R}^m$ \hbx
       \end{rema}

       Since $\pi_1M$  is a finite abelian group, it can be represented as
       $\pi_1M= G(q_1)\oplus \cdots \oplus G(q_s)$ where $G(q) = \mathbb{Z}_{q^{\alpha_1}} \oplus \cdots \oplus \mathbb{Z}_{q^{\alpha_u}}$ , $\alpha_i\in \mathbb{N}$ and  where $q_1 <\cdots < q_u$ are prime numbers, see   (\ref{DecompositionOfpi_1(M)}).

   Since $f_{\#}=identity$,
  $\mathcal{R}(f^k)= \pi_1M/im({\rm id -id}) =\pi_1M= G(q_1)\oplus \cdots \oplus G(q_u)$
  for all $k|n$. We notice that $L(f^k)=0$ $\iff$ $d_i|k$ for some $i=1,...,u$.

 Moreover the action of $f_{\#}={\rm id}$ on $\mathcal{R}(f^k)$ is constant, hence each orbit contains a single class and $\mathcal{OR}(f^k)=\mathcal{R}(f^k)$.

   To prove the remaining Lemmas it will be convenient to formulate the assumptions which are satisfied for self-maps of compact connected Lie groups $f:M\to M$ satisfying ${\rm dim}M\ge 3$ and all eigenvalues of $A(f)$ are roots of unity.

\begin{assum}

 We consider a self-map $f:M\to M$ of a compact    rational exterior manifold  of dimension $m\ge 3$ satisfying

 \begin{enumerate}
 \item  There exists a matrix $A\in \mathcal{M}_{r\times r}(\mathbb{Z})$  whose all eigenvalues have modulus $= 1$,  and  $L(f^n)=det(I-A^n)$ for all $n\in \mathbb{N}$.
  \item    $\pi_1M$ is finite commutative  $\pi_1M= G(q_1)\oplus \cdots \oplus G(q_s)$
  and the homotopy homomorphism $f_{\#}= identity$

  \item $f$ is a Jiang map i.e. for each $k$ all classes in $\mathcal{R}(f^k)$ have the same index.
\hbx
  \end{enumerate}
\end{assum}
\begin{rema}
The Assumptions (1)-(3) are satisfied by a self-map $f:M\to M$ of a compact Lie group satisfying the right hand side of Theorem \ref{MainTeoIdentity} (i.e. $dim M\ge 3$ , all eigenvalues of $A(f)$ have modulus $\le 1$  and $f_{\#}=id$).  In fact  $A$  may be     the matrix of the homomorphism $A(f)/\ker f$ see Remark \ref{Less=1Implies=1}

\end{rema}

  \begin{lema}\label{d=p^alpha}
     The Assumptions (1)-(3) imply  the following.
     If $L(f)\ne 0$ then  for each $j=1,...,u$ there exist $i=1,...,v$ and $\beta\in \mathbb{N}$ so that $d_i=q_j^{\beta}$

  \end{lema}
  {\it Proof}. Since $q_j$ is a prime, it is enough to show that for each $j=1,...,u$ there exist $i=1,...,v$ and $\beta\in \mathbb{N}$ so that $d_i$ divides $q_j^{\beta}$.

  Assume contrary i.e. there exists $j$ so that $  d_i$ does not divide $q_j^{\beta}$ for all
    $i=1,...,v$ and $\beta \in \mathbb{N}$. Now (1) implies $L(f^{q_j^{\beta}})\ne 0$ for all $\beta \in \mathbb{N}$,      hence by
    the Jiang property  all classes in $\mathcal{R}(f^{q_j^k})$ are
    essential.
    On the other hand  the homomorphism $i_{q_j^k,q_j^l}:\mathcal{R}(f^{q_j^l})
    \to \mathcal{R}(f^{q_j^k})$ is given by
    $$ [x]\longrightarrow [(id+f_{\#}^{q_j^l}+f_{\#}^{2\cdot q_j^l}+\cdots +
    f_{\#}^{q_j^k-q_j^l})(x)] = [x+\cdots +x]=[q_j^k/q_j^l\cdot x]=[q_j^{k-l}\cdot x]$$
    Here we used the assumption that $f_{\#}=id$.

    This implies
    $im i_{q_j^k,q_j^l}\subset im i_{q_j^k,q_j^{k-1}}\varsubsetneq \mathcal{R}(f^{q_j^k})$ (for $l<k)$, hence there are irreducible  classes. Index of each such class
is divisible by $q_j^k$. By the Jiang property $L(f^{q_j^k})={\rm ind}(f^{q_j^k})$ is also divisible by $q_j^k$.
  Since by our assumption no $d_i$ divides $q_j^k$, $L(f^{q_j^k})\ne 0$ which implies $\lim_{k \to \infty}|L(f^{q_j^k})|=+\infty$.
  The last gives a contradiction, since $L(f^k)$ is cyclic. \hbx

  \begin{lema}\label{i_kl=bijection}
   Under  the above Assumptions (1) and (2)

  \centerline{$i_{kl}:\mathcal{R}(f^l) \to \mathcal{R}(f^k)$  is a bijection $\iff$
  no  $q_j$  divides $k/l$.}

  \end{lema}

{\it Proof}. Let us recall that  $\pi_1 M= G(q_1)\oplus \cdots \oplus G(q_u)$ and the map $i_{kl}: \mathcal{R}(f^l) \to \mathcal{R}(f^k)$
$$i_{kl}[x]=(id+f^l_{\#}+f^{2\cdot l}_{\#}+\cdots + +f^{k-l}_{\#})[x]
=[x+\cdots+x]=[k/l\cdot x]$$
 splits to the product of self-maps of
all summands, each map given by the same formula as above.
It remains to check when the restriction od $i_{kl}$ to each summand is the bijection. But  the homomorphism
$\mathbb{Z}_{q^{\alpha}}\ni x \to hx\in \mathbb{Z}_{q^{\alpha}}$
is a bijection exactly  for $q\nmid h$.\hbx

We will denote $\bar k={\rm gcd}(k,d)$.

\begin{lema}\label{bar_kIZO}
 Under the Assumptions (1) and (2).
If $L(f^k)\ne 0$ then $i_{k\bar k}$ is an isomorphism.

\end{lema}

{\it Proof}. By Lemma \ref{i_kl=bijection} it remains to show that no $q_j$ divides $k/\bar k$.

We fix a $j=1,...,s$. By Lemma \ref{d=p^alpha},   $d_i=q_j^{\beta}$  for an $i=1,...,v$; $\beta \in \mathbb{N}$.
Let $d=q_j^{\beta_0}\cdot d'$ where $q_j\nmid d'$ and
let $k=q_j^{\beta_1}\cdot k'$ where $q_j \nmid k'$.

Since  $d_i=q_j^{\beta}$ divides $d=q_j^{\beta_0}\cdot d'$,
$\beta\le \beta_0$.
On the other hand  $L(f^k)\ne 0$ so $d_i=q_j^{\beta}$ does not divide $k=q_j^{\beta_1}\cdot k'$, hence
 $\beta_1<\beta$. This implies  $\beta_1<\beta_0$ and
 $$\bar k = gcd(k,d)=gcd(q_j^{\beta_1}\cdot k', q_j^{\beta_0}\cdot d')=
 q_j^{\beta_1}\cdot gcd (k',  d')$$
 Now $$k/\bar k=(q_j^{\beta_1}\cdot k')
 /(q_j^{\beta_1}\cdot gcd (k',  d'))= k'
 /  gcd (k',  d')$$
 Since $q_j$ does not divide $k'$, it neither divides $k/\bar k$.
 \hbx

 \begin{coro}\label{a|D}
  Under the  Assumptions (1) , (2) ,and (3).
 If $A\in \mathcal{R}(f^a)$ is an essential irreducible class then $a|d$

 \end{coro}

 {\it Proof}. $A\in \mathcal{IEOR}$ and $f^â$ is Jiang imply  $L(f^a)\ne 0$. Now, by Lemma
 \ref{bar_kIZO}, $i_{a,\bar a}$ is a bijection, hence $A$ reduces to a (unique) class in $\mathcal{R}(f^{\bar a})$ where $\bar a=lcm(a,d)$. Since $A$ is irreducible, $a=\bar a$ and the last divides $d$.\hbx




    \begin{lema}\label{L=L^bar}Under the Assumption (1).
    $L(f^k)=L(f^{\bar k})$ for $\bar k={\rm gcd}(k,d)$.

    \end{lema}

    {\it Proof}. This follows from Lemma 4.3 in \cite{Sinica} where it is shown that under the above assumptions the Dold expansion is finite
    $L(f^k)=\sum_{l|d}a_l\cdot {\rm reg}_l(k)$. This implies
    $$L(f^k)=\sum_{l|d}a_l\cdot {\rm reg}_l(k) = \sum_{l|{\rm gcd}(k,d)}a_l\cdot {\rm reg}_{l}(k) =\sum_{l|\bar k} a_l\cdot l  = L(f^{\bar k})$$  \hbx

   \begin{exem}{\rm
    We give an example illustrating the above theory. Let us consider the special unitary group $SU(2)$. Its dimension is $3$ and it may be identified with the group of quaternions of the unit length. Its center consists of two elements $\{+1; -1\}$.
The quotient  ${\rm PSU}(2)=  {\rm SU}(2)/\{+1,-1\}$ is a Lie group and can be identified with the real projective space $\mathbb{RP}^3$, hence  ${\rm PSU}(2)$ is an orientable manifold and $\pi_1(PSU(2))=\mathbb{Z}_2$. Moreover  $H^q({\rm PSU}(2);\mathbb{Q})= H^q(\mathbb{RP}^3;\mathbb{Q})=H^q(S^3;\mathbb{Q}) = \mathbb{Q}$ for $q=0,3$ and zero otherwise. Now  $H^q({\rm PSU}(2);\mathbb{Q})=\Lambda^*(a)$ is the exterior algebra  and  $0\ne a\in H^3(\rm PSU(2);\mathbb{Q})$.

Let $M=\rm PSU(2)\times \rm PSU(2)\times \dots \times \rm PSU(2)$ (s times). Then $H^*(M;\mathbb{Q})= \Lambda^*(a_1,...,a_s; \mathbb{Q})$ which implies $A(M)=\mathbb{Q}(a_1,...,a_s)$ where the last means the space with the  basis
$(a_1,...,a_s)$. Moreover $\pi_1M= \mathbb{Z}_2 \times \cdots \mathbb{Z}_2$.

Let $f:M\to M$ be a continuous map. For fixed numbers $1\le k,l\le s$ we define  inclusions $i_l: {\rm PSU(2)} \to M$ on l-th component $i_l(x)=(1,...,1,x,1,...,1)$ and projections $p_k:M\to  {\rm PSU(2)}$ given by $p_k(x_1,...,x_s)= x_k$. Let $\alpha_{lk}= {\rm deg}(p_k f i_l)$ be the degree of the composition $p_k f i_l: {\rm PSU(2)} \to  {\rm PSU(2)}$ . Now
$$f^*(a_k)=f^*p_k^*(a)= (i_1^*f^*p_k^*(a), i_2^*f^*p_k^*(a), ... , i_s^*f^*p_k^*(a))$$
$$=  (p_kfi_1)^*(a), (p_kfi_2)^*(a), ... , (p_kfi_s)^*(a)=(\alpha_{1k}\cdot a_1,...,\alpha_{sk}\cdot a_s)$$
so $A(f): A(M) \to A(M)$ is given by the matrix $[\alpha_{lk}]$.

On the other hand we notice that for each integer matrix $[\alpha_{lk}]$ there exists a map $f:M\to M$ such that $p_kfi_l= \alpha_{lk}$. In fact we define the maps
$f_{lk}: {\rm PSU}(2)\to {\rm PSU}(2)$ of degree $\alpha_{lk}$ and then $f_k: M ={\rm PSU}(2) \times \cdots \times {\rm PSU}(2) \to {\rm PSU}(2)$ by
$f_k(x_1,...,x_s)=(f_{1k}(x_1)\cdot \cdots \cdot (f_{sk}(x_s)))$ (here $\cdot$ denotes the multiplication. Finally we put $f:M\to M$ as
$f(x)= (f_1(x),...,f_s(x))$).}
\end{exem}

   \end{document}